\documentclass[12pt,reqno,a4wide]{amsart}

\usepackage{tikz} 
\usepackage{calrsfs}
\usepackage[charter]{mathdesign}
\usepackage{hyperref}

\allowdisplaybreaks

\begin{document}

\title{Sums of squares of Tetranacci numbers: A generating function approach}

\author[Helmut Prodinger]{Helmut Prodinger}
\author[Sarah J. Selkirk]{Sarah J. Selkirk}
\address{Department of Mathematics, University of Stellenbosch 7602,
	Stellenbosch, South Africa}
\email{hproding@sun.ac.za}
\email{sjselkirk@sun.ac.za}

\begin{abstract} 
	It is demonstrated how an explicit expression of the (partial) sum of Tetranacci numbers can be found and proved using generating functions and the Hadamard product. We also provide a Binet-type formula for generalized Fibonacci numbers, by explicitly factoring the denominator of their generating functions.
\end{abstract}

\maketitle

\section{Introduction}

\emph{Tetranacci numbers} $u_n$ (OEIS: A000078, \cite{OEIS}) are defined either by the recursion 
\begin{align*}
u_{n+4}=u_{n+3}+u_{n+2}+u_{n+1}+u_{n}, \qquad \text{where }\qquad u_0=0,\ u_1=1,\ u_2=1,\ u_3=2,
\end{align*}
or via the generating function
\begin{equation*}
\sum_{n\ge0}u_nz^n=\frac{z}{1-z-z^2-z^3-z^4}.
	\end{equation*} 

A typical result in the recent paper \cite{schumacher} is the evaluation
\begin{equation*}
\sum_{0\le k\le n}u_k^2=\frac13+u_nu_{n+1}
-\frac13(u_{n+1}-u_{n-1})^2+\frac13u_nu_{n-2} +\frac13u_{n-2}u_{n-3}.
\end{equation*}
for which a (long) proof by induction had been given.

The present note wants to shed some light on how to use generating functions to prove such a result and also how to find this (or an equivalent formula).

Furthermore, \emph{all} the roots of the polynomial $1 - z - z^2 - \cdots - z^h$ are explicitly determined in terms of generalized binomial series. This leads to a Binet-type formula for generalized Fibonacci numbers, and as an example we provide the formula for Tetranacci numbers. 

\section{The Hadamard product of two power series}

For two power series (generating functions) $f(z)=\sum_n a_nz^n$, $g(z)=\sum_n b_nz^n$, the Hadamard product is defined as
\begin{equation*}
\sum_{n\ge0}a_nb_nz^n.
\end{equation*}
If both $f(z)$ and $g(z)$ are rational, the resulting power series of their Hadamard product is again rational. There are computer algorithms to do this effectively, for instance \textsc{Gfun} \cite{GFUN}, implemented in \textsc{Maple}. 

We first provide a simple example of the Hadamard product of two generating functions: Let
\begin{equation*}
f(z)=\sum_{n\ge0}n2^nz^n=\frac{2z}{(1-2z)^2},\qquad \text{and}\qquad g(z)=\sum_{n\ge0}n3^nz^n=\frac{3z}{(1-3z)^2}.
\end{equation*}
By multiplying corresponding coefficients of $f(z)$ and $g(z)$, the Hadamard product is given by 
\begin{equation*}
\sum_{n\ge0}n^26^nz^n=\frac{6z(1+6z)}{(1-6z)^2},
\end{equation*}
where the generating function is computed by \textsc{Gfun}.

To give an example in the context of Tetranacci numbers, the Hadamard product of the generating function $z/(1-z-z^2-z^3-z^4)$ with itself is given by
\begin{equation*}
\sum_{n\ge0}u_n^2z^n=\frac{z-{z}^{2}-2\,{z}^{3}-2\,{z}^{4}-2\,{z}^{5}+{z}^{6}+{z}^{7}}
{1-2\,z-4\,{z}^{2}-6\,{z}^{3}-12\,{z}^{4}+4\,{z}^{5}+6\,{z}^{6}+2\,{z}
	^{8}-{z}^{10}
			}.
	\end{equation*}
By general principles, the generating function of the partial sums is then given as
\begin{align}\label{l1}
\sum_{n\ge0}\Big(\sum_{0\le k\le n}u_k^2\Big)z^n&=\frac1{1-z}\cdot\frac{z-{z}^{2}-2\,{z}^{3}-2\,{z}^{4}-2\,{z}^{5}+{z}^{6}+{z}^{7}}
{1-2\,z-4\,{z}^{2}-6\,{z}^{3}-12\,{z}^{4}+4\,{z}^{5}+6\,{z}^{6}+2\,{z}
	^{8}-{z}^{10}
}
\nonumber\\
&=\frac{z}{3(1-z)}+\frac{2\,z+{z}^{2}-{z}^{3}-{z}^{4}+5\,{z}^{5}+4\,{z}^{6}+{z}^{7}+{z}^{8}-{z
	}^{9}-{z}^{10}
	}{3(1-2\,z-4\,{z}^{2}-6\,{z}^{3}-12\,{z}^{4}+4\,{z}^{5}+6\,{z}^{6}+2\,{z}
	^{8}-{z}^{10})
	}
\end{align}
For a simpler expression, we can compute the Hadamard product involving coefficients $u_n$ and $u_{n+2}$,
\begin{equation*}
\sum_{n\ge0}u_nu_{n+2}z^n=\frac{2z}
{1-2\,z-4\,{z}^{2}-6\,{z}^{3}-12\,{z}^{4}+4\,{z}^{5}+6\,{z}^{6}+2\,{z}
	^{8}-{z}^{10}},
\end{equation*}
and by letting $t_n=\frac12u_nu_{n+2}$, we can compare this with the generating function in (\ref{l1}) to find that for $n\ge1$,
\begin{equation*}
	\sum_{0\le k\le n}u_k^2=\frac13+\frac13(2t_{n}+t_{n-1}-t_{n-2}-t_{n-3}+5t_{n-4}+4t_{n-5}+t_{n-6}+t_{n-7}-t_{n-8}-t_{n-9}).
\end{equation*}
This is an equivalent formula for the one obtained in \cite{schumacher}. Note that $t_n=0$ for negative indices.
The fact these two formulas are indeed equivalent can be checked by a computer, and all the generating functions
\begin{equation*}
\sum_{n\ge0}u_{n-i}u_{n-j}z^n,
\end{equation*}
with fixed integers $i$, $j$,
can be effectively computed via the Hadamard product algorithm implemented in \textsc{Gfun}.

\section{Higher order Fibonacci-type recursions}

To show how the generating function machinery works on similar but more involved sums, let us step up a bit and define
\begin{equation*}
\sum_{n\ge0}u_nz^n=\frac{z}{1-z-z^2-z^3-z^4-z^5}.
\end{equation*}
Again computing the Hadamard product of this generating function with itself, we find
\begin{align*}
	&\sum_{n\ge0}\Big(\sum_{0\le k\le n}u_k^2\Big)z^n=
	\frac{3z}{8(1-z)}\\&+\frac{\scriptstyle{-5\,z-3\,{z}^{2}+{z}^{3}+4\,{z}^{4}+2\,{z}^{5}-34\,{z}^{6}-30\,{z}^{7
		}-20\,{z}^{8}-20\,{z}^{9}-16\,{z}^{10}+6\,{z}^{11}+6\,{z}^{12}+3\,{z}^
		{13}+3\,{z}^{14}+3\,{z}^{15}}
						}{8(1-2 z-4 {z}^{2}-7 {z}^{3}-14 {z}^{4}-28 {z}^{5}+4 {z}^{6}+6 {z
		}^{7}+4 {z}^{9}+10 {z}^{10}-{z}^{12}-{z}^{15}
		)}
	\end{align*}
With a shift in coefficients as done in the previous example, we compute that
\begin{equation*}
\sum_{n\ge0}u_nu_{n+3}z^n=\frac{4z}{1-2 z-4 {z}^{2}-7 {z}^{3}-14 {z}^{4}-28 {z}^{5}+4 {z}^{6}+6 {z
}^{7}+4 {z}^{9}+10 {z}^{10}-{z}^{12}-{z}^{15}}.
\end{equation*}
We let $t_n=\frac{1}{4}u_nu_{n+3}$ and can then express the sum in question as follows:
\begin{equation*}
\sum_{0\le k\le n}u_k^2=\frac38+\frac18(-5t_n-3t_{n-1}+t_{n-2}+4t_{n-3}+\cdots+3t_{n-13}+3t_{n-14}).
\end{equation*}

This process can be generalized to any higher order Fibonacci-type recursion, such as 
\begin{equation*}
	\sum_{n\ge0}u_nz^n=\frac{z}{1-z-z^2-z^3-z^4-z^5-z^6},
\end{equation*}
and other identities and related expressions can also be computed, but are too long to be displayed here.

\section{Higher order Fibonacci-type numbers}

For generalized Fibonacci numbers defined by the usual initial values and the recursion
\begin{align*}
u_{n+h} = u_{n+h-1} + u_{n+h-2} + \cdots + u_{n}, 
\end{align*}
the corresponding generating function (and its simplification) is
\begin{equation*}
\frac{z}{1-z-\cdots-z^h}=\frac{z}{1-z\cdot\dfrac{1-z^h}{1-z}}=\frac{z(1-z)}{1-2z+z^{h+1}}.
\end{equation*}

The dominant root of this rational function already occurs in the literature, see for example \cite{Prodinger}. However, we can do better than that and describe \emph{all} the roots of the denominator, obtaining in this
way a Binet-type formula. We consider the generating function
\begin{equation*}
\frac{1}{1-2z+z^{h+1}},
\end{equation*}
from which the original case can be obtained by simple shifts. 

We determine the roots of the denominator in terms of generalized binomial series, going back to Lambert, and described in more detail in \cite{GKP}. A generalized binomial series is defined as  
\begin{align*}
\big(\mathcal{B}_{t}(x)\big)^r = \sum_{n \geq 0}\binom{tn+r}{n}\frac{r}{tn+r}x^n.
\end{align*}
Given the expression $1-\frac{z}{u}+z^{h+1}$, let $\zeta$ be a primitive $h$-th root of unity. Then the $h+1$ roots can be expressed in terms of these generalized binomial series as 
\begin{align*}
u\mathcal{B}_{h+1}(u^{h+1}) \qquad \text{and} \qquad \zeta^{-j}u^{-\frac{1}{h}}\mathcal{B}_{(h+1)/h}\big(\zeta^j u^{\frac{h+1}{h}}\big)^{-\frac{1}{h}} \quad \text{for $0 \leq j \leq h-1$}.
\end{align*}
In our case, we are dealing with the special case where $u = \frac{1}{2}$. 

It is easy to verify (and the calculation for $u = \frac{1}{2}$ has appeared in \cite{Prodinger}) that 
\begin{align*}
1 - u\mathcal{B}_{h+1}(u^{h+1}) + u^{h+1}\mathcal{B}_{h+1}(u^{h+1})^{h+1} & = 0.
\end{align*}
Now, the other roots can be checked by considering first $j = 0$:
\begin{align*}
&1 - u^{-(h+1)/h}\mathcal{B}_{(h+1)/h}\big(u^{\frac{h+1}{h}}\big)^{-\frac{1}{h}} + u^{-\frac{(h+1)}{h}}\mathcal{B}_{(h+1)/h}\big(u^{\frac{h+1}{h}}\big)^{-\frac{(h+1)}{h}}\\
& = 1 - \sum_{n \geq 0}\binom{\frac{(h+1)n - 1}{h}}{n}\frac{-1}{(h+1)n - 1}u^{\frac{(h+1)(n-1)}{h}} + \sum_{n \geq 0}\binom{\frac{(h+1)(n-1)}{h}}{n}\frac{-1}{n - 1}u^{\frac{(h+1)(n-1)}{h}}.
\end{align*}
Since 
\begin{align*}
\binom{\frac{(h+1)n - 1}{h}}{n}\frac{-1}{(h+1)n - 1} & = -\frac{(\frac{(h+1)(n-1)}{h})\cdots (\frac{n - 1 + h}{h})}{h \cdot n!},
\end{align*}
and 
\begin{align*}
\binom{\frac{(h+1)(n-1)}{h}}{n}\frac{-1}{n - 1} & = -\frac{(\frac{(h+1)(n-1)}{h})\cdot (\frac{(h+1)(n-1) - h}{h}) \cdots (\frac{n - 1 + h}{h}) \cdot (\frac{n-1}{h})}{(n-1)n!} = -\frac{(\frac{(h+1)(n-1)}{h})\cdots (\frac{n - 1 + h}{h})}{h \cdot n!},
\end{align*}
the result follows. The roots for $j \neq 0$ follow from the substitution $u = 1\cdot u$, and with the power of $\frac{1}{h}$ playing a role at each $u$, we obtain all possible $h$-th roots of unity.  

From this, we can explicitly compute the coefficients of 
\begin{align*}
\frac{1}{1-2z+z^{h+1}}.
\end{align*}
For ease of notation, let $r_h = \frac{1}{2}\mathcal{B}_{h+1}\Big(\frac{1}{2^{h+1}}\Big)$, and for $0 \leq j \leq h-1$ let
\begin{align*}
r_j = \zeta^{-j}2^{\frac{1}{h}}\mathcal{B}_{(h+1)/h}\Big(\zeta^j \Big(\frac{1}{2}\Big)^{\frac{h+1}{h}}\Big)^{-\frac{1}{h}}.
\end{align*}
Then using partial fractions and these $r_i$ values, we can compute that $[z^n]\frac{1}{1-2z+z^{h+1}}$ is equal to 
\begin{align*}
[z^n]\frac{1}{(z-r_0)(z-r_1)\cdots(z-r_h)} = [z^n]\sum_{i = 0}^h\frac{1}{(z-r_i)}\prod_{\substack{j = 0\\ j \neq i}}^{h}(r_i - r_j)^{-1} = -\sum_{i = 0}^h \frac{1}{r_i^{n+1}}\prod_{\substack{j = 0\\ j \neq i}}^{h}(r_i - r_j)^{-1}.
\end{align*}
Therefore we have obtained a Binet-type formula for generalized Fibonacci numbers. In fact, with these coefficients it is possible to compute and verify identities on the level of coefficients for expressions such as those discussed in previous sections.   


\subsection{A formula for Tetranacci numbers}

To provide a concrete example of how one would use these roots to compute generalized Fibonacci numbers, we provide the calculation of the formula for the case with the Tetranacci numbers. 

Tetranacci numbers correspond to $h = 4$, so the five roots are (as calculated by a computer):
\begin{align*}
r_4 & = \frac{1}{2}\mathcal{B}_{5}\Big(\frac{1}{2^{5}}\Big) = 0.518790063675884,\\
r_0 & = 2^{\frac{1}{4}}\mathcal{B}_{5/4}\Big(\Big(\frac{1}{2}\Big)^{\frac{5}{4}}\Big)^{-\frac{1}{4}} = 1,\\
r_1 & = -i2^{\frac{1}{4}}\mathcal{B}_{5/4}\Big(i \Big(\frac{1}{2}\Big)^{\frac{5}{4}}\Big)^{-\frac{1}{4}} = -0.114070631164587 -1.21674600397435i,\\
r_2 & = -2^{\frac{1}{4}}\mathcal{B}_{5/4}\Big(- \Big(\frac{1}{2}\Big)^{\frac{5}{4}}\Big)^{-\frac{1}{4}} = -1.29064880134671,\\
r_3 & = i2^{\frac{1}{4}}\mathcal{B}_{5/4}\Big(-i \Big(\frac{1}{2}\Big)^{\frac{5}{4}}\Big)^{-\frac{1}{4}} = -0.114070631164587 + 1.21674600397435i.
\end{align*}
Using these roots and the initial values, we can determine the values of $A$, $B$, $C$, $D$, and $E$ in the expression
\begin{align*}
u_n = A\cdot r_4^{-n} + B\cdot r_0^{-n} + C\cdot r_1^{-n} + D\cdot r_2^{-n} + E\cdot r_3^{-n}.
\end{align*}
Again using a computer, we find that these are given by
\begin{align*}
A & = 0.293813062773642\\
B & = 0\\
C & = -0.0504502052166080 - 0.169681902881564i\\
D & = -0.192912652340427\\
E & = -0.0504502052166080 + 0.169681902881564i
\end{align*}
Therefore the $n$-th Tetranacci number can be calculated via the formula:
\begin{align*}
u_n &= \frac{0.293813062773642}{(0.518790063675884)^n} - \frac{0.0504502052166080 + 0.169681902881564i}{(-0.114070631164587 -1.21674600397435i)^n}\\
&\quad - \frac{0.192912652340427}{(-1.29064880134671)^n} + \frac{-0.0504502052166080 + 0.169681902881564i}{(-0.114070631164587 + 1.21674600397435i)^n}.
\end{align*}
Analogous computations provide similar formulas for other generalized Fibonacci numbers. 

\newpage

\bibliographystyle{plain}

\end{document}